\documentclass[a4paper,10pt]{article}
\usepackage{amsmath}
\usepackage{amssymb}
\usepackage{amsthm}
\theoremstyle{plain}
\numberwithin{equation}{section}
\newtheorem{theo}{Theorem}[section]
\newtheorem{prop}{Proposition}[section]
\newtheorem{cor}{Corollary}[section]
\newtheorem{lem}{Lemma}[section]
\theoremstyle{definition}
\newtheorem{defi}{Definition}[section]
\theoremstyle{remark}
\newtheorem{rem}{Remark}[section]
\newtheorem{exm}{Example}[section]
\begin{document}
\begin{Large}
 {\bf Second order subexponentiality and infinite\\ divisibility}
\end{Large}
\begin{center}
 Toshiro Watanabe\\
  
\medskip
{\bf Abstract}
  
\medskip

\end{center}  
 We characterize  the second order subexponentiality of an infinitely divisible distribution on the real  line under an exponential moment assumption. We investigate the asymptotic behaviour of the difference between the tails of an infinitely divisible distribution and  its L\'evy measure. Moreover, we  study  the second order asymptotic behaviour of  the tail of the $t$-th convolution power of an infinitely divisible distribution. The density version for a self-decomposable distribution on the real  line
without an exponential moment assumption is also  given. Finally, the regularly varying case for a self-decomposable distribution on the half  line is discussed.
  

subexponentiality,
  infinite divisibility \\

\section{Introduction and results}
\medskip
The subexponentiality of infinitely divisible distributions on the half line was characterized by Embrechts et al.\ \cite{egv} and on the real line by Pakes \cite{p}.  The subexponentiality of an infinitely divisible distribution implies the asymptotic equivalence between the tails of the distribution and  its L\'evy measure. In this paper, we characterize the second order subexponentiality of an infinitely divisible distribution on the real  line in terms of its L\'evy measure under an exponential moment assumption. The  second order subexponentiality  yields a higher asymptotic relation than the usual subexponentiality between the tails of an infinitely divisible distribution and  its L\'evy measure. 

In what follows, we denote by $\mathbb R$ the real line and by $\mathbb R_{+}$ the half line $[0,\infty)$. Denote by $\mathbb N$ the totality of positive integers. The symbol $\delta_a(dx)$ stands for the delta measure at $a \in \mathbb R$. Let $\eta$ and $\rho$ be probability distributions on $\mathbb R$. We denote by $\eta*\rho$ the convolution of $\eta$ and $\rho$ and by $\rho^{n*}$ the
$n$-th convolution power of $\rho$ with the understanding that $\rho^{0*}(dx)=\delta_0(dx)$. Denote by $m(\rho)$ the mean of $\rho$. The characteristic function of $\rho$ is denoted by $\widehat\rho(z)$, namely, for $z \in \mathbb R$, 
\begin{equation}
\widehat\rho(z):=\int_{-\infty}^{\infty}e^{izx}\rho(dx). \nonumber
\end{equation}
 For a measure $\xi$ on $\mathbb R$, we denote by $\bar\xi(x)$ the tail $\xi((x,\infty))$ for $x >0$. For positive functions $f(x)$ and $g(x)$ on $[a,\infty)$ for some $a \in \mathbb R$, we define the relation $f(x)\sim g(x)$ by $\lim_{x \to \infty}f(x)/g(x)=1.$ We say that $f(x,A)\sim c  f(x)$ as $x \to \infty$ and then $A \to \infty$,  if
  $$\lim_{A \to \infty}\lim_{x \to \infty}f(x,A)/f(x)=c>0.  $$
 We say that  $f(x,A)=o(f( x))$ as $x \to \infty$ and then $A \to \infty$, if
   $$\limsup_{A \to \infty}\limsup_{x \to \infty}|f(x,A)|/f(x)=0.   $$    
\begin{defi} (i)  A nonnegative measurable function $g(x)$ on $\mathbb R$ belongs to the class ${\bf L}$ if $g(x+a)\sim g(x)$ for every $ a \in \mathbb R$.

 (ii)  Let $\Delta:=(0,c]$ with $c >0$. A distribution $\rho$ on $\mathbb R$ belongs to the class $\mathcal{L}_{\Delta}$ if $\rho((x,x+c]) \in {\bf L}$. A distribution $\rho$ on $\mathbb R$ belongs to the class $\mathcal{L}_{loc}$ if $\rho \in \mathcal{L}_{\Delta}$ for each $\Delta:=(0,c]$ with $c >0$.

(iii)  Let $\Delta:=(0,c]$ with $c >0$. A distribution $\rho$ on $\mathbb R$ belongs to the class $\mathcal{S}_{\Delta}$ if $\rho \in \mathcal{L}_{\Delta}$ and  $\rho^{2*}((x,x+c])\sim 2\rho((x,x+c])$. 
 A distribution $\rho$ on $\mathbb R$ belongs to the class $\mathcal{S}_{loc}$ if $\rho \in \mathcal{S}_{\Delta}$ for each $\Delta:=(0,c]$ with $c >0$.
\end{defi}
If a distribution $\rho$ on $\mathbb R$ belongs to the class $\mathcal{L}_{loc}$, then,  for  $c >0$, 
$$\rho((x,x+c])\sim c\rho((x,x+1])$$
and, for every $\delta >0$, $e^{\delta x}\rho((x,x+1]) \to\infty$ as $x \to \infty$. See  (2.6) in the proof of Theorem 2.1  of Watanabe and Yamamuro \cite{wy2} and Lemma 2.17 of Foss et al.\ \cite{fkz}.
A distribution $\rho$ on $\mathbb R$ belongs to the class $\mathcal{S}$ if $\bar\rho(x) \in {\bf L}$ and $\overline{\rho^{2*}}(x) \sim 2\overline{\rho}(x).$ Distributions in the classes
$\mathcal{S}$ and  $\mathcal{S}_{loc}$
are called {\it subexponential} and  {\it locally
 subexponential}, respectively. 
\medskip
\begin{defi} A distribution $\rho$ on $\mathbb R$ belongs to the class $\mathcal{S}^2_{loc}$ if the following three conditions hold :

(1) $\rho \in \mathcal{S}_{loc}$.  

 (2) $\int_{-\infty}^{\infty}|x|\rho(dx) < \infty$.

(3) We have
\begin{equation}
\overline{\rho^{2*}}(x)=2\bar\rho(x)+2m(\rho)\rho((x,x+1])+o(\rho((x,x+1]))
\end{equation}
 as $x \to \infty.$

\end{defi}
The subclasses $\mathcal{S}_{\Delta}$, $\mathcal{S}_{loc}$, and $ \mathcal{S}^2_{loc}$  of the class $\mathcal{S}$ were respectively  introduced by Asmussen et al.\ \cite{afk}, Watanabe and Yamamuro \cite{wy2}, and  Lin \cite{l}.  Lin \cite{l} treated the one-sided case and used the
 symbol $\mathcal{S}_2$ for the class $ \mathcal{S}^2_{loc}$. Distributions in the class $\mathcal{S}^2_{loc}$ are called {\it second order subexponential}. Infinitely divisible distributions
 on  $\mathbb R$
  in the classes $\mathcal{S}_{\Delta}$
and
$\mathcal{S}_{loc}$ are
found in Watanabe and Yamamuro  \cite{wy1,wy2} and Shimura and Watanabe \cite{sw}. Lin \cite{l} gave some sufficient conditions in order that a distribution  on $\mathbb R_+$ belongs to the class $\mathcal{S}^2_{loc}$. See Proposition 2.4 and Corollary 2.1 of \cite{l}. He showed that the  lognormal distribution,  Weibull distribution with parameter $\beta \in(0,1)$, and Pareto distribution with parameter $\alpha >1$ belong to the class $\mathcal{S}^2_{loc}$.
 Geluk and  Pakes \cite{gp} and  Geluk \cite{g} treated another second order subexponentiality.
\medskip
 Let $\mu$ be an infinitely divisible distribution on $\mathbb R$. Then, its characteristic function $\widehat\mu(z)$ is represented as
\begin{equation}
\widehat\mu(z)=\exp \left( \int_{-\infty}^{\infty}(e^{izx}-1-\frac{izx}{1+x^2})\nu(dx)+i\gamma z-\frac{1}{2}az^2 \right),         
\nonumber
\end{equation}
where $\gamma \in \mathbb R$, $a \geq 0$, and $\nu$ is a measure on $\mathbb R$ satisfying $\nu(\{0\})=0$ and 
$$\int_{-\infty}^{\infty}\frac{x^2}{1+x^2}\nu(dx) < \infty.$$
The measure $\nu$ is called L\'evy measure of $\mu$. See Sato  
 \cite{s}. Throughout the paper, we assume that  the tail $\bar\nu(c) $ is positive for all $c > 0$. For $c >0$, define a normalized distribution $\nu_{(c)}$ as
 $$ \nu_{(c)}(dx):=1_{(c,\infty)}(x)\frac{\nu(dx)}{\bar\nu(c)}.$$
Here the symbol $1_{(c,\infty)}(x)$ stands for the indicator function of the  set $(c,\infty)$. Denote by $\mu^{t*}$ the $t$-th convolution power of $\mu$ for $t >0$. Note that $\mu^{t*}$ is the distribution of $X_t$ for a certain L\'evy process $\{X_t\}$. 
\begin{theo}
Let $\mu$ be an  infinitely divisible distribution on $\mathbb R$  with L\'evy measure $\nu$. Assume that there exists $\epsilon >0$ such that $\int_{-\infty}^{\infty}\exp(-\epsilon x)\mu(dx) < \infty$.   Then, we have the  following : 

(i) $\mu \in \mathcal{S}^2_{loc}$ if and only if $\nu_{(1)}  \in \mathcal{S}^2_{loc}$.

(ii) If $\mu \in \mathcal{S}^2_{loc}$, then
\begin{equation}
\bar\nu(x)=\bar\mu(x)-m(\mu)\mu((x,x+1])+o(\mu((x,x+1]))
\end{equation} 
as $x \to \infty,$ equivalently,
\begin{equation}
\bar\mu(x)=\bar\nu(x)+m(\mu)\nu((x,x+1])+o(\nu((x,x+1]))
\end{equation} 
as $x \to \infty.$

(iii) Conversely, if (1.2) with finite $m(\mu)$, $\mu \in \mathcal{S}_{loc}$, and 
$$(\bar\mu(x))^2=o(\mu((x,x+1]))$$
 as $x \to \infty$ hold, then $\mu \in \mathcal{S}^2_{loc}$.
\end{theo}
\begin{rem} An exponential moment assumption  in the above theorem is  necessary for the restriction of the class  $\mathcal{S}_{loc}$ in the two sided case. See Jian et al.\ \cite{jwcc}
for the detailed account.

\end{rem}
\begin{cor}
Let $\mu$ be an  infinitely divisible distribution on $\mathbb R$  with L\'evy measure $\nu$. Assume that there exists $\epsilon >0$ such that $\int_{-\infty}^{\infty}\exp(-\epsilon x)\mu(dx) < \infty$.   Then, we have the following :  

(i) $\mu \in \mathcal{S}^2_{loc}$ if and only if $\mu^{t*} \in \mathcal{S}^2_{loc}$ for some $t>0$, equivalently, for all $t>0$.

(ii) If $\mu \in \mathcal{S}^2_{loc}$, then, for all $t>0$,
\begin{equation}
\overline{\mu^{t*}}(x)=t\bar\mu(x)+(t^2-t)m(\mu)\mu((x,x+1])+o(\mu((x,x+1]))
\end{equation} 
as $x \to \infty.$ 
\end{cor}
\begin{rem} Let $\mu$ be an  infinitely divisible distribution on $\mathbb R_+$  with L\'evy measure $\nu$. If $\mu \in \mathcal{S}_{loc}$, $m(\mu) < \infty$, and $\mu$ satisfies (1.4) for $t=t_0, t_0+1$ with some $t_0>0$, then  $\mu \in \mathcal{S}^2_{loc}$. 
\end{rem}

The organization of this paper is as follows. In Sect.\ 2, we give  preliminaries for the proof of the results. In Sect.\ 3, we prove Theorem 1.1 and its corollary. In Sect.\ 4, we treat the self-decomposable case. In Sect.\ 5, three examples of the results are given.
In Sect.\ 6, we give some remarks on the regularly varying case.
\section{Preliminaries}

Watanabe and Yamamuro \cite{wy2} used the main results of Watanabe \cite{w} on the convolution equivalence of infinitely divisible distributions on $\mathbb R$  to prove the following two lemmas.

\begin{lem} (Corollary 2.1 of \cite{wy2}) Let $\mu$ be an  infinitely divisible distribution on $\mathbb R$  with L\'evy measure $\nu$. Assume that there exists $\epsilon >0$ such that $\int_{-\infty}^{\infty}\exp(-\epsilon x)\mu(dx) < \infty$. Then, the following are equivalent : 

(1) $\mu \in \mathcal{S}_{loc}$.

(2) $\nu_{(1)} \in \mathcal{S}_{loc}$.

(3) $\nu_{(1)} \in \mathcal{L}_{loc}$
and $\mu((x,x+c]) \sim \nu((x,x+c])$ for all $c >0$.

\end{lem}
\begin{rem} Since $\mathcal{S}_{loc} \subset \mathcal{S}$, we see that if condition (1) holds in the above lemma, then 
$$\bar\mu(x)\sim \bar\nu(x) \in {\bf L}.$$

\end{rem}
\begin{lem} (Corollary 3.1 of \cite{wy2}) Let $\mu$ be an  infinitely divisible distribution on $\mathbb R$  with L\'evy measure $\nu$.  Assume that there exists $\epsilon >0$ such that $\int_{-\infty}^{\infty}\exp(-\epsilon x)\mu(dx) < \infty$.
 If $\mu^{t*} \in \mathcal{S}_{loc}$ for some $t >0$, then $\mu^{t*} \in \mathcal{S}_{loc}$ for all $t >0$ and
$$ \mu^{t*}((x,x+c]) \sim t\mu((x,x+c]) $$
 for all $t >0 $ and for all $c >0. $

\end{lem}
Lin \cite{l} proved the following three lemmas.
\begin{lem}(Theorem 2.1 of  \cite{l})
Let $\rho$ be a distribution on $\mathbb R_+$. Let $\{p_n\}_{n=0}^{\infty}$ be a nonnegative sequence with $p_n >0$ for some $n \geq 2$ and $\sum_{n=0}^{\infty}p_n =1$ satisfying $\sum_{n=0}^{\infty}p_n(1+\epsilon_1)^n < \infty$ for some $\epsilon_1 >0$.
Define a distribution $\eta$ on $\mathbb R_+$ as
$$ \eta(dx):=\sum_{n=0}^{\infty}p_n\rho^{n*}(dx).$$
Then we have the following:

(i) If $\rho \in \mathcal{S}^2_{loc},$ then we have $\eta \in \mathcal{S}^2_{loc}$ and
\begin{equation}
\overline{\eta}(x)=(\sum_{n=1}^{\infty}np_n)\bar\rho(x)+(\sum_{n=2}^{\infty}n(n-1)p_n)m(\rho)\rho((x,x+1])+o(\rho((x,x+1]))
\end{equation}
 as $x \to \infty.$

(ii) Conversely, if (2.1) with finite $m(\rho)$, $\rho \in \mathcal{S}_{loc},$ and 
$$(\bar\rho(x))^2=o(\rho((x,x+1])$$
 as $x \to \infty$ hold, then  $\rho \in \mathcal{S}^2_{loc}.$
\end{lem}
\begin{rem} We can see from the proof of Theorem 2.1 of  \cite{l} that even in the case of $p_n <0$ for some $ n \geq 0$, assertion (i) of  the above lemma is still true if $\sum_{n=0}^{\infty}|p_n|(1+\epsilon_1)^n < \infty$ for some $\epsilon_1 >0$.

\end{rem}
\begin{lem}(Proposition 2.3 of  \cite {l})
Let $\rho$ and  $\eta$ be
 distributions on $\mathbb R_+$. If $\rho \in \mathcal{S}^2_{loc},$ and there are $K>0$ and $c \in \mathbb R$ such that
 $$\lim_{x \to \infty}\frac{\bar\eta(x)-K\bar\rho(x)}{\rho((x,x+1])}=c,$$
 then $\eta \in \mathcal{S}^2_{loc}.$
\end{lem}

\begin{lem}(Lemma 3.4 of  \cite{l}) Let $\rho$  be a
 distribution on $\mathbb R_+$. Assume that $m(\rho) < \infty$, $\rho \in \mathcal{L}_{loc}$ and $(\bar\rho(x))^2=o(\rho(x,x+1])$ as $x \to \infty$. Then the relation (1.1)
implies $\rho \in \mathcal{S}_{loc}$. 
\end{lem}
Let $\delta:=\bar\nu(c)$ for $c>0$. Define a compound Poisson distribution $\mu_1$ and a distribution $\sigma$ on $\mathbb R_+$ as
\begin{equation}
\mu_1:=e^{-\delta}\sum_{n=0}^{\infty}\frac{\delta^n}{n!}(\nu_{(c)})^{n*}
\end{equation}
and
\begin{equation}
\sigma:=\frac{e^{-\delta}}{1-e^{-\delta}}\sum_{n=1}^{\infty}\frac{\delta^n}{n!}(\nu_{(c)})^{n*}.
\end{equation}
\begin{lem}
We can choose sufficiently large $c>0$ such that $0 <e^{\delta}-1<1$ and we have
\begin{equation}
\nu_{(c)}=-\frac{1}{\delta}\sum_{n=1}^{\infty}\frac{(1-e^{\delta})^n}{n}\sigma^{n*}.
\end{equation}
\end{lem}
Proof. We define a a signed measure  $\eta$ as
\begin{equation}
\eta:=-\frac{1}{\delta}\sum_{n=1}^{\infty}\frac{(1-e^{\delta})^n}{n}\sigma^{n*}. \nonumber
\end{equation}
Let $\rho$ be a signed measure on $\mathbb R_+$.
Denote by $L_{\rho}(t)$ for $t \geq 0$ the Laplace transform of $\rho$, that is, $L_{\rho}(t):=\int_{0-}^{\infty}e^{-tx}\rho(dx).$
We have
\begin{equation}
\begin{split}
L_{\eta}(t)&=-\frac{1}{\delta}\sum_{n=1}^{\infty}\frac{(1-e^{\delta})^n}{n}(L_{\sigma}(t))^n\\
&=\frac{1}{\delta}\log(1-(1-e^{\delta})L_{\sigma}(t)).\nonumber
\end{split}
\end{equation}
We see from (2.3) that
$$L_{\sigma}(t)=(e^{\delta}-1)^{-1}(\exp(\delta L_{\nu_{(c)}}(t))-1).$$
Thus we have
\begin{equation}
\begin{split}
L_{\eta}(t)&=\frac{1}{\delta}\log(\exp(\delta L_{\nu_{(c)}}(t)))\\
&=L_{\nu_{(c)}}(t)\nonumber
\end{split}
\end{equation}
and hence we have $\eta=\nu_{(c)}$,
that is, (2.4).
 \hfill $\Box$
\section{Proof of Theorem 1.1 and its corollary}
Proof of Theorem 1.1. Let $\mu$ be an infinitely divisible distribution on $\mathbb R$ with L\'evy measure $\nu$. As in Lemma 2.6, we choose sufficiently large $c>0$ such that $0 <e^{\delta}-1<1$. We define an infinitely divisible distribution $\mu_2$ by $\mu=\mu_1*\mu_2$.    Assume that there exists $\epsilon \in (0,1)$ such that $\int_{-\infty}^{\infty}\exp(-\epsilon x)\mu(dx) < \infty$.  Then we see from Theorem 25.17  of Sato \cite{s} that $\int_{-\infty}^{\infty}\exp(-\epsilon x)\mu_2(dx) < \infty$ and, for every $b>0$,  $\int_{-\infty}^{\infty}\exp(b x)\mu_2(dx) < \infty$. Hence, for every $b>0$,  $\overline{\mu_2}(x) = o(e^{-bx})$ as $x \to \infty$. We find from Lemma 2.1 that  $\mu \in \mathcal{S}_{loc}$  if and only if
$\mu_1 \in \mathcal{S}_{loc}.$  Since $\int_{-\infty}^{\infty}\exp(-\epsilon x)\mu(dx) < \infty$, we have $\int_{-\infty}^{\infty}|x|\mu(dx) < \infty$ if and only if $\int_{-\infty}^{\infty}|x|\mu_1(dx) < \infty$. Suppose that $\mu \in \mathcal{S}_{loc}$ and $\int_{-\infty}^{\infty}|x|\mu(dx) < \infty$, that is, $\mu_1 \in \mathcal{S}_{loc}$ and $\int_{-\infty}^{\infty}|x|\mu_1(dx) < \infty$. We have
\begin{equation}
\begin{split}
&\bar\mu(x)-\bar\mu_1(x)\\
&=\overline{\mu_1*\mu_2}(x)-\bar\mu_1(x)\\
&=\int_{0-}^{\infty}\mu_1((x-y,x])\mu_2(dy)-\int_{-\infty}^{0-}\mu_1((x,x-y])\mu_2(dy)\\
&=I_1-I_2.
\end{split}
\end{equation}
If $\int_{0-}^{\infty}y\mu_2(dy)=0$, then $I_1=0$ and if 
$\int_{-\infty}^{0- } |y |\mu_2(dy)=0$, then $I_2=0$. Thus we can assume that $\int_{0-}^{\infty}y\mu_2(dy)>0$ and $\int_{-\infty}^{0- } |y |\mu_2(dy)>0$.
We find that
$$I_1:=I_{11}+I_{12}+I_{13}.$$
where, for $A >0$, 
$$I_{11}:=\int_{0-}^{A+}\mu_1((x-y,x])\mu_2(dy),$$
$$I_{12}:=\int_{A+}^{ x/2+}\mu_1((x-y,x])\mu_2(dy),$$
and
$$I_{13}:=\int_{ x/2+}^{\infty}\mu_1((x-y,x])\mu_2(dy).$$ 
We have by $\mu_1 \in  \mathcal{S}_{loc}\subset \mathcal{L}_{loc}$
\begin{equation}
\begin{split}
I_{11}&\sim \mu_1((x,x+1])\int_{0-}^{A+}y\mu_2(dy)\\
&\sim \mu_1((x,x+1])\int_{0-}^{\infty}y\mu_2(dy) \nonumber
\end{split}
\end{equation}
 as $x \to \infty$ and then  $A \to \infty$. For any $\epsilon_1 \in (0,1)$, there is $C_1>0$ such that, for $0 \leq y \leq  x/2$ and for sufficiently large $x >0$,
$$\frac{\mu_1((x-y,x])}{\mu_1((x,x+1])}\leq C_1e^{\epsilon_1 y}.$$
Thus we see that
\begin{equation}
\begin{split}
I_{12}&\leq \mu_1((x,x+1]) C_1\int_{A+}^{ x/2+}e^{\epsilon_1 y}\mu_2(dy)\\
&=o(\mu_1((x,x+1])),\nonumber
\end{split}
\end{equation}
 as $x \to \infty$ and then  $A \to \infty$.
We have  
$$I_{13}\leq\bar\mu_2( x/2)=o(e^{-x})=o(\mu_1((x,x+1]))$$
 as $x \to \infty$. Thus we have
 \begin{equation}
I_{1}\sim \mu_1((x,x+1])\int_{0- }^{\infty} y \mu_2(dy).
\end{equation}For any $\epsilon_2 \in (0,\epsilon)$, there is $C_2>0$ such that, for $ y< 0$ and for sufficiently large $x >0$,
$$\frac{\mu_1(x,,x-y])}{\mu_1((x,x+1])}\leq C_2e^{\epsilon_2 |y|}.$$
Thus, by dominated convergence theorem, we see that
\begin{equation}
I_{2}\sim \mu_1((x,x+1])\int_{-\infty}^{0- } |y |\mu_2(dy).
\end{equation}
Hence, we find from (3.1), (3.2), and (3.3) that 
\begin{equation}
\bar\mu(x)-\bar\mu_1(x)=m(\mu_2)\mu_1((x,x+1])+o(\mu_1((x,x+1]))
\end{equation}
as $x \to \infty$.
By argument analogous to the above equation, we have
$$\overline{\mu^{2*}}(x)-\overline{\mu_1^{2*}}(x)=m(\mu_2^{2*})\mu_1^{2*}((x,x+1])+o(\mu_1^{2*}((x,x+1]))$$
as $x \to \infty$.
Since $\mu_1^{2*}((x,x+1])\sim 2\mu_1((x,x+1])$ and $m(\mu_2^{2*})=2m(\mu_2)$, 
$$\overline{\mu^{2*}}(x)-\overline{\mu_1^{2*}}(x)=4m
(\mu_2)\mu_1((x,x+1])+o(\mu_1((x,x+1]))$$
as $x \to \infty$.
Thus we see from (3.4) that
\begin{equation}
\begin{split}
&\overline{\mu^{2*}}(x)-2\overline{\mu}(x)\\
&=\overline{\mu_1^{2*}}(x)-2\overline{\mu_1}(x)+2m(\mu_2)\mu_1((x,x+1])+o(\mu_1((x,x+1]))
\end{split}
\end{equation}
 as $x \to \infty$. Since $m(\mu)=m(\mu_1)+m(\mu_2)$ and we find from Lemma 2.1 that
\begin{equation}
\mu((x,x+1]))\sim\nu((x,x+1]))\sim\mu_1((x,x+1]),
\end{equation}
we have by (3.5)
$$\overline{\mu^{2*}}(x)=2\bar\mu(x)+2m(\mu)\mu((x,x+1])+o(\mu((x,x+1]))$$
 as $x \to \infty$ if and only if
 $$\overline{\mu_1^{2*}}(x)=2\bar\mu_1(x)+2m(\mu_1)\mu_1((x,x+1])+o(\mu_1((x,x+1]))$$
 as $x \to \infty$. Hence $\mu \in \mathcal{S}^2_{loc}$ if and only if $\mu_1  \in \mathcal{S}^2_{loc}$. Since, for $x > 0$, 
 $$\frac{\overline{\mu_1}(x)}{1-e^{-\delta}}=\bar\sigma(x)$$
we see from Lemma 2.4 that $\sigma \in \mathcal{S}^2_{loc}$ if and only if $\mu_1  \in \mathcal{S}^2_{loc}$.
We find from (2.3), Lemma 2.6, and Remark 2.2 that if $\sigma \in \mathcal{S}^2_{loc}$, then $\nu_{(c)} \in \mathcal{S}^2_{loc}$ for sufficiently large $c >0$. We see from (2.2) and Lemma 2.3 that if  $\nu_{(c)} \in \mathcal{S}^2_{loc}$, then $\mu_1  \in \mathcal{S}^2_{loc}$. Thus, for sufficiently large $c >0$, 
$\mu \in \mathcal{S}^2_{loc}$ if and only if $\nu_{(c)} \in \mathcal{S}^2_{loc}$. Since, for sufficiently large $x > 0$,
$$\overline{\nu_{(c)}}(x) =\frac{\bar\nu(1)}{\bar\nu(c)}\overline{\nu_{(1)}}(x),$$
we obtain from Lemma 2.4 that $\nu_{(1)} \in \mathcal{S}^2_{loc}$ if and only if $\nu_{(c)} \in \mathcal{S}^2_{loc}$ for sufficiently large $c >0$. Thus we have $\mu \in \mathcal{S}^2_{loc}$ if and only if $\nu_{(1)} \in \mathcal{S}^2_{loc}$. We have proved assertion (i).
Next, we prove assertion (ii). Assume that $\mu \in \mathcal{S}^2_{loc}$, equivalently, $\nu_{(c)} \in \mathcal{S}^2_{loc}$ for $c >0$. Note that $m(\mu_1)=\delta m(\nu_{(c)})$. We see from Lemma 2.3 that
\begin{equation}
\begin{split}
\overline{\mu_1}(x)&=e^{-\delta}\sum_{n=1}^{\infty}\frac{\delta^n}{(n-1)!}\overline{\nu_{(c)}}(x) \\
&+e^{-\delta}\sum_{n=2}^{\infty}\frac{\delta^n}{(n-2)!}m(\nu_{(c)})\nu_{(c)}((x,x+1])+o(\nu_{(c)}((x,x+1]))\\
&=\overline{\nu}(x)+m(\mu_1)\nu((x,x+1])+o(\nu((x,x+1]))
\end{split}
\end{equation}
as $x \to \infty.$
Thus we obtain (1.2) and (1.3) from (3.4) and (3.6). Next we prove assertion (iii). We see from (3.4) that the assumption that
(1.2) with finite $m(\mu)$, $\mu \in \mathcal{S}_{loc}$, and $(\bar\mu(x))^2=o(\mu((x,x+1]))$ as $x \to \infty$
is equivalent to that (3.7) with finite $m(\nu_{(c)})$, $\nu_{(c)} \in \mathcal{S}_{loc}$, and $(\overline{\nu_{(c)}}(x))^2=o(\nu_{(c)}((x,x+1]))$ as $x \to \infty$. This implies from Lemma 2.3 that $\nu_{(c)} \in \mathcal{S}^2_{loc}$, equivalently, $\mu \in \mathcal{S}^2_{loc}$. \hfill $\Box$

\medskip
Proof of Corollary 1.1. We see from Theorem 1.1 that 
 $\mu^{t*}
  \in \mathcal{S}^2_{loc}$ for some $t >0$,   equivalently, for all $t >0$ if   and only if $\nu_{(1)}  \in \mathcal{S}^2_{loc}$.
Hence assertion (i) is true. Next we prove assertion (ii). Suppose that  $\mu
  \in \mathcal{S}^2_{loc}$. Then we find from (i) that
  $\mu^{t*}
  \in \mathcal{S}^2_{loc}$ for all $t >0$. We see from (1.2) that
$$\overline{\mu^{t*}}(x)=t\bar\nu(x)+m(\mu^{t*})\mu^{t*}((x,x+1])+o(\mu^{t*}((x,x+1]))$$
 as $x \to \infty$.
 Note that $m(\mu^{t*})=tm(\mu)$ and from Lemma 2.2 that 
$$\mu^{t*}((x,x+1])\sim t\mu((x,x+1]).$$
 Thus we have  by (1.2)
\begin{equation}
\begin{split}
\overline{\mu^{t*}}(x)&=t\bar\nu(x)+t^2m(\mu)\mu((x,x+1])+o(\mu((x,x+1]))\\
&=t\bar\mu(x)+(t^2-t)m(\mu)\mu((x,x+1])+o(\mu((x,x+1]))\nonumber
\end{split}
\end{equation}
 as $x \to \infty$.
 We have proved (1.4). \hfill $\Box$

Proof of Remark 1.2. Assume that $\mu \in \mathcal{S}_{loc}$, $m(\mu) < \infty$, and $\mu$ satisfies (1.4) for $t=t_0, t_0+1$ with some $t_0>0$. Then we have
 \begin{equation}
\begin{split}
&\overline{\mu^{(t_0+1)*}}(x)-(t_0+1)\overline{\mu}(x)
-t_0(\overline{\mu^{2*}}(x)-2\overline{\mu}(x)) \\
&=\int_{0-}^{x+}\overline{\mu^{t_0*}}(x-y)\mu(dy)+\overline{\mu}(x)
-(t_0+1)\overline{\mu}(x) \\
&-t_0\int_{0-}^{x+}\overline{\mu}(x-y)\mu(dy)-t_0\overline{\mu}(x)
+2t_0\overline{\mu}(x) \\
&=\int_{0-}^{x+}(\overline{\mu^{t_0*}}(x-y)-t_0\overline{\mu}(x-y))
\mu(dy)\\
&= I_1 +I_2 +I_3,\\
\end{split}
\end{equation}
where, for $0 <2 A <x$,
$$ I_1:=\int_{0-}^{A+}(\overline{\mu^{t_0*}}(x-y)-t_0\overline{\mu}(x-y))
\mu(dy),$$
$$ I_2:=\int_{A+}^{(x-A)+}(\overline{\mu^{t_0*}}(x-y)-t_0\overline{\mu}(x-y))
\mu(dy),$$
and
$$ I_3:=\int_{(x-A)+}^{x+}(\overline{\mu^{t_0*}}(x-y)-t_0\overline{\mu}(x-y))
\mu(dy).$$
We divide the proof into three cases: $t_0>1$; $t_0=1$; and $0<t_0<1$.
Let $t_0>1$. By the assumption, we see that
\begin{equation}
\begin{split}
I_1&\sim t_0(t_0-1)m(\mu)\int_{0-}^{A+}\mu((x-y,x-y+1])\mu(dy)\\
&\sim   t_0(t_0-1)m(\mu)\mu((x,x+1])
\end{split}
\end{equation}
as $x \to \infty$ and then $A \to \infty$.
We find from  $\mu \in \mathcal{S}_{loc}$ that there is $\epsilon>0$
such that
\begin{equation}
\begin{split}
|I_2|&\leq (1+\epsilon)t_0(t_0-1)m(\mu)\int_{A+}^{(x-A)+}\mu((x-y,x-y+1])\mu(dy)\\
&=  o(\mu((x,x+1]))
\end{split}
\end{equation}
as $x \to \infty$ and then  $A \to \infty$.
By using integration by parts, we have
\begin{equation}
\begin{split}
 I_3&=\int_{0-}^{A+}(\bar\mu(x-y)-\bar\mu(x))\mu^{t_0*}(dy)\\
& - t_0\int_{0-}^{A+}(\bar\mu(x-y)-\bar\mu(x))\mu(dy)\\
&+(\overline{\mu^{t_0*}}(A)- t_0\bar\mu(A))(\bar\mu(x-A)-\bar\mu(x))\\
&=K_1 -K_2 +K_3. \nonumber
\end{split}
\end{equation}
As  $x \to \infty$ and then $A \to \infty$, we have
$$K_1\sim m(\mu^{t_0*})\mu((x,x+1])= t_0 m(\mu)\mu((x,x+1]),$$
and
$$K_2\sim  t_0 m(\mu)\mu((x,x+1]).$$
Note from $m(\mu)<\infty$
 that  $\overline{\mu^{t_0*}}(A)A \to 0$ and $\overline{\mu}(A)A \to 0$ as   $A \to \infty$.
Thus we see that 
$$\limsup_{A \to \infty}\limsup_{x \to \infty}\frac{|K_3|}{\mu((x,x+1])}\leq \limsup_{A \to \infty}(\overline{\mu^{t_0*}}(A)+ t_0\bar\mu(A))A=0.$$
 Thus we have
\begin{equation}
 I_3=o(\mu((x,x+1]))
\end{equation}
as $x \to \infty$ and then $A \to \infty$. Thus we obtain from (3.8), (3.9), (3.10), and (3.11) and the assumption that
\begin{equation}
\begin{split}
&(t_0+1)t_0m(\mu)\mu((x,x+1])-t_0(\overline{\mu^{2*}}(x)-2\overline{\mu}(x))\\
& = t_0(t_0-1)m(\mu)\mu((x,x+1]) +o(\mu((x,x+1]))
\end{split}
\end{equation}
 as $x \to \infty$. Hence we have
\begin{equation}
\overline{\mu^{2*}}(x)=2\overline{\mu}(x)+2m(\mu)\mu((x,x+1])+o(\mu((x,x+1])
\end{equation}
as $x \to \infty$. That is,  $\mu \in \mathcal{S}^2_{loc}$. 
Next, let $t_0=1$. Then we have (3.13) and hence $\mu \in \mathcal{S}^2_{loc}$. Finally, let $0<t_0<1$. In the same way, we see that, as 
$x \to \infty$ and then $A \to \infty$,
$$-I_1\sim t_0(1-t_0)m(\mu)\mu((x,x+1]),$$
$$I_2=  o(\mu((x,x+1])),$$
and
$$I_3=  o(\mu((x,x+1])).$$
Thus we have (3.12) and (3.13) and hence $\mu \in \mathcal{S}^2_{loc}$.
\hfill $\Box$

\section{Self-decomposable case}

Let $f(x)$ and $g(x)$ be probability density functions on $\mathbb R$. We denote by $f \otimes g(x)$ the convolution of $f(x)$ and $g(x)$ and by $f^{n\otimes}(x)$ the $n$-th convolution power of $f(x)$ for $n \in \mathbb N$.
\begin{defi} (i) A probability density function $g(x)$ on $\mathbb R$ belongs to the class $\mathcal{L}_d$ if $g(x) \in {\bf L}$. 

(ii)  A probability density function $g(x)$ on $\mathbb R$ belongs to the class $\mathcal{S}_d$ if $g(x) \in \mathcal{L}_d$ and $g^{2\otimes}(x)\sim 2g(x)$. 
\end{defi}

\begin{defi} A  probability density function $g(x)$ on $\mathbb R$ belongs to the class $\mathcal{S}^2_{d}$ if the following three conditions hold :

(1) $g(x) \in \mathcal{S}_{d}$.  

 (2) $\int_{-\infty}^{\infty}|x|g(x)dx < \infty$.

(3) For $\rho(dx):=g(x)dx$,
$$\overline{\rho^{2*}}(x)=2\bar\rho(x)+2m(\rho)g(x)+o(g(x)),$$
 as $x \to \infty.$

\end{defi}
The classes $\mathcal{S}_d$ and $\mathcal{S}^2_{d}$ were introduced by Chover et al.\ \cite{cnw} and Omey and Willekens \cite{ow1}, respectively. Densities in the classes $\mathcal{S}_d$ and $\mathcal{S}^2_{d}$ are called {\it  subexponential} and {\it second order subexponential}, respectively. 
See also Foss et al.\ \cite{fkz} and Kl\"uppelberg \cite{k2}
 for the class $\mathcal{S}_d$.  An infinitely divisible distribution on  $\mathbb R_+$ with its density in the class  $\mathcal{S}_d$ is found in Watanabe \cite{w3}. Omey and Willekens \cite{ow1} studied an infinitely divisible distribution on  $\mathbb R_+$ with the density of the normalized L\'evy measure in the class $\mathcal{S}^2_{d}$. However, they could not characterize the density of an infinitely divisible distribution on  $\mathbb R_+$ with its density in the class  $\mathcal{S}^2_d$ because they did not know Lemma 2.1 and Lemma 4.1 below.
An infinitely divisible distribution $\mu$ on $\mathbb R$ is called {\it self-decomposable} if, for every $b \in (0,1)$, there is a distribution $\rho_b$ on $\mathbb R$ such that 
$$\widehat\mu(z)=\widehat\mu(bz)\widehat{\rho_b}(z).$$
An infinitely divisible distribution $\mu$ on $\mathbb R$ is self-decomposable if and only if $\nu(dx)=k(x)/|x|dx$ with $k(x)$ being nonnegative and increasing on $(-\infty,0)$ and nonnegative and decreasing on $(0,\infty)$. An infinitely divisible distribution $\mu$ on $\mathbb R$ is non-degenerate if it is not a delta measure. Every non-degenerate self-decomposable  distribution $\mu$ on $\mathbb R$ is absolutely continuous and unimodal. Many important statistical distributions are known to be self-decomposable. However their L\'evy measures and the $t$-th convolution powers are often not explicitly known. See Sato \cite{s}.
Let $\mu(dx)=p(x)dx$ be a non-degenerate self-decomposable  distribution on $\mathbb R$. We assume that $k(x)$ is positive for all $x >0$. We define self-decomposable distributions $\xi_1(dx)=p_1(x)dx
$ and $\xi_2(dx)=p_2(x)dx$ as $\mu=\xi_1*\xi_2$ and
\begin{equation}
\widehat\xi_1(z):=\exp \left( \int_{0}^{\infty}(e^{izx}-1)\frac{k(x\vee d)}{x}dx\right).         
\nonumber
\end{equation}
for sufficiently large $d>0$. Watanabe and Yamamuro \cite{wy2} proved the following two lemmas.
\begin{lem}(Theorem 1.3 of \cite{wy2} and its proof) Let $\mu(dx)=p(x)dx$ be a self-decomposable  distribution on $\mathbb R$ with $\nu(dx)=k(x)/|x|dx$. The following are equivalent : 

(1)  $\mu \in \mathcal{S}_{loc}$.

(2) $p(x) \in \mathcal{S}_{d}$.

(3) $p_1(x) \in \mathcal{S}_{d}$.

(4) $\frac{1}{\bar\nu(1)}1_{(1,\infty)}(x)k(x)/x \in \mathcal{S}_{d}$.

(5) $k(x) \in {\bf L}$ and $p(x)\sim p_1(x) \sim k(x)/x$.
\end{lem}

\begin{rem} Let $\mu(dx)=p(x)dx$ be a self-decomposable  distribution on $\mathbb R$ with L\'evy measure $\nu(dx)=k(x)/|x|dx$.
We see from Lemma 4.1 that $\mu \in \mathcal{S}^2_{loc}$ if and only if  $p(x) \in \mathcal{S}^2_{d}$ and that $\nu_{(1)}  \in \mathcal{S}^2_{loc}$ if and only if $\frac{1}{\bar\nu(1)}1_{(1,\infty)}(x)k(x)/x \in \mathcal{S}^2_{d}$.
\end{rem}
\begin{lem} (Theorem 1.4 of \cite{wy2}) Let  $\mu(dx)=p(x)dx$ be a self-decomposable  distribution on $\mathbb R$. Let $p^t(x)$ be the density of $\mu^{t*}$ for  $t >0$.
 If $p^t(x) \in \mathcal{S}_{d}$ for some $t >0$, then $p^t(x) \in \mathcal{S}_{d}$ for all $t >0$ and
$$ p^t(x) \sim tp(x) $$
 for all $t >0. $
\end{lem}
\begin{prop} Let $\mu(dx)=p(x)dx$ be a self-decomposable  distribution on $\mathbb R_+$. If  $(\bar\mu(x))^2 =o(\mu((x,x+1]))$ as $x \to \infty$, $m(\mu) < \infty$, and
$$\overline{\mu^{2*}}(x)=2\bar\mu(x)+2m(\mu)\mu((x,x+1])+o(\mu((x,x+1]))$$
 as $x \to \infty,$ then $p(x) \in \mathcal{S}^2_{d}$. 
\end{prop}
Proof. Assume that  $(\bar\mu(x))^2 =o(\mu((x,x+1]))$ as $x \to \infty$, $m(\mu) < \infty$, and
$$\overline{\mu^{2*}}(x)=2\bar\mu(x)+2m(\mu)\mu((x,x+1])+o(\mu((x,x+1]))$$
 as $x \to \infty.$  Note that
$$\overline{\mu^{2*}}(x)-2\bar\mu(x)+(\bar\mu(x))^2=\int_{0-}^{x+}(\bar\mu(x-y)-\bar\mu(x))\mu(dy).$$
Thus, by the assumption, we have
$$\int_{0-}^{x+}(\bar\mu(x-y)-\bar\mu(x))\mu(dy)\sim 2m(\mu)\mu((x,x+1]).$$
We shall prove that, for every $m \in \mathbb N,$
\begin{equation}
\lim_{x \to \infty}\frac{\mu((x-m,x-m+1])}{\mu((x,x+1])}=1.
\end{equation}
 Since $\mu$ is unimodal, we see that, for every $m \in \mathbb N,$
 $$\liminf_{x \to \infty}\frac{\mu((x-m,x-m+1])}{\mu((x,x+1])}\geq 1.$$
  Suppose that there are some $c >1$, $m_0 \in \mathbb N$, and a increasing sequence $\{x_n\}_{n=1}^{\infty}$ with $\lim_{n \to \infty}x_n= \infty$ such that
 $$\lim_{n \to \infty}\frac{\mu((x_n-m_0,x_n-m_0+1])}{\mu((x_n,x_n+1])}=c.$$
we have
$$\int_{0-}^{x_n+}(\bar\mu(x_n-y)-\bar\mu(x_n))\mu(dy)=I_1+I_2
+I_3,
$$
where
$$I_1:=\int_{0-}^{m+}(\bar\mu(x_n-y)-\bar\mu(x_n))\mu(dy),$$
$$I_2:=\int_{m+}^{(x_n-m)+}(\bar\mu(x_n-y)-\bar\mu(x_n))\mu(dy),$$
and
$$I_3:=\int_{(x_n-m)+}^{x_n+}(\bar\mu(x_n-y)-\bar\mu(x_n))\mu(dy).$$
By the unimodality, we have for $0 \leq y \leq m_0$
 $$\liminf_{n \to \infty}\frac{\mu((x_n-y,x_n])}{\mu((x_n,x_n+1])}\geq  y$$
and  for $y \geq m_0$
 $$\liminf_{n \to \infty}\frac{\mu((x_n-y,x_n])}{\mu((x_n,x_n+1])}\geq  c(y-m_0)+m_0$$
Thus we have by Fatou's lemma
\begin{equation}
\begin{split}
&\liminf_{m \to \infty}\liminf_{n \to \infty}\frac{I_1}{\mu((x_n,x_n+1])}\\
&\geq\int_{0-}^{m_0+}\liminf_{n \to \infty}\frac{\mu((x_n-y,x_n])}{\mu((x_n,x_n+1])}\mu(dy)\\
&+\liminf_{m \to \infty}\int_{m_0+}^{m+}\liminf_{n \to \infty}\frac{\mu((x_n-y,x_n])}{\mu((x_n,x_n+1])}\mu(dy)\\
&\geq\int_{0-}^{m_0+}y\mu(dy)+
 \liminf_{m \to \infty}\int_{m_0+}^{m+}( c(y-m_0)+m_0)\mu(dy)\\
&= \int_{0-}^{m_0+}y\mu(dy)+
 \int_{m_0+}^{\infty}( c(y-m_0)+m_0)\mu(dy)
 >m(\mu).
 \end{split}
\end{equation}
Clearly, we have
\begin{equation}
\liminf_{n \to \infty}\frac{I_2}{\mu((x_n,x_n+1])}\geq 0.
\end{equation}
By using integration by parts, we see that, for sufficiently large $n$,
\begin{equation}
I_3=I_1+(\bar\mu(x_n-m)-\bar\mu(x_n))(\bar\mu(m)-\bar\mu(x_n))\geq I_1. \nonumber
\end{equation}
Thus we obtain from (4.2) that
\begin{equation}
\liminf_{m \to \infty}\liminf_{n \to \infty}\frac{I_3}{\mu((x_n,x_n+1])}>m(\mu).
\end{equation}
Hence we have by (4.2), (4.3), and (4.4)
$$\liminf_{n \to \infty}\frac{1}{\mu((x_n,x_n+1])}\int_{0-}^{x_n+}(\bar\mu(x_n-y)-\bar\mu(x_n))\mu(dy)>2m(\mu).$$
This is a contradiction. Thus we have proved (4.1).
By the unimodality, it implies that $p(x) \in  \mathcal{L}_{d}$. Thus, by Lemma 2.5, we have proved
$\mu \in  \mathcal{S}^2_{loc}$ and hence $p(x) \in  \mathcal{S}^2_{d}$. \hfill $\Box$
\begin{theo}
Let $\mu(dx)=p(x)dx$  be a self-decomposable  distribution on $\mathbb R$  with L\'evy measure $\nu(dx)=k(x)/|x|dx$. Assume that $\int_{-\infty}^{0-}|x|\mu(dx) < \infty$.   Then, the following hold : 

(i) $p(x) \in \mathcal{S}^2_{d}$ if and only if $\nu_{(1)}  \in \mathcal{S}^2_{loc}$, equivalently, $\frac{1}{\bar\nu(1)}1_{(1,\infty)}(x)k(x)/x \in \mathcal{S}^2_{d}$.

(ii) If $p(x) \in \mathcal{S}^2_{d} $, then
\begin{equation}
\bar\nu(x)=\bar\mu(x)-m(\mu)p(x)+o(p(x))
\end{equation} 
as $x \to \infty,$ equivalently,
\begin{equation}
\bar\mu(x)=\bar\nu(x)+m(\mu)k(x)/x+o(k(x)/x)
\end{equation} 
as $x \to \infty.$

(iii) Conversely, if (4.5) with finite $m(\mu)$, $p(x) \in \mathcal{S}_{d} $, and $(\bar\mu(x))^2=o(p(x))$ as $x \to \infty$ hold, then $p(x) \in \mathcal{S}^2_{d} $.
\end{theo}
Proof.  Since the support of the L\'evy measure of $\xi_2$ has an  upper bound, we find from Theorem 25.17 of Sato \cite{s} that,
for every $b >0$, $\int_0^{\infty}e^{bx}\xi_2(dx) < \infty$ and hence $\overline{\xi_2}(x) = o(e^{-bx})$ as $x \to \infty$.
We have 
\begin{equation}
\begin{split}
&\bar\mu(x)-\bar\xi_1(x)\\
&=\overline{\xi_1*\xi_2}(x)-\bar\xi_1(x)\\
&=\int_{0}^{\infty}\xi_1((x-y,x])\xi_2(dy)-\int_{-\infty}^{0}\xi_1((x,x-y])\xi_2(dy)\\
&=I_1-I_2.\nonumber
\end{split}
\end{equation}
If $\int_{0}^{\infty}y\xi_2(dy)=0$, then $I_1=0$ and if 
$\int_{-\infty}^{0 } |y |\xi_2(dy)=0$, then $I_2=0$. Thus we can assume that $\int_{0}^{\infty}y\xi_2(dy)>0$ and $\int_{-\infty}^{0 } |y |\xi_2(dy)>0$.
Suppose that $p(x) \in \mathcal{S}_{d}$, equivalently by Lemma 4.1, $p_1(x) \in \mathcal{S}_{d}$.  Note that $p_1(x) \in {\bf L}$ and $\xi_1$ is unimodal. Thus there are $C>0$ and $\epsilon >0$ such that, for $0 <y < x/2$ and for sufficiently large $x >0$,
\begin{equation}
\xi_1((x-y,x]) \leq Ce^{\epsilon y}p_1(x).\nonumber
\end{equation}
Note that $\int_{0}^{\infty}e^{\epsilon y}\xi_2(dy) < \infty$ and
$$\int_{x/2}^{\infty}\xi_1((x-y,x])\xi_2(dy)\leq \bar\xi_2(x/2)=o(e^{-x})=o(p_1(x))$$
as $x \to \infty.$ Thus, by dominated convergence theorem,
$$I_1 \sim p_1(x)\int_{0}^{\infty}y\xi_2(dy).$$
Since $\xi_1$ is unimodal, we have, for $y <0$ and for sufficiently large $x >0$,
$$\xi_1((x,x-y])\leq p_1(x)|y|.$$
Since $\int_{-\infty}^{0}|y|\mu(dy)< \infty$, we see from Theorem 25.3 of Sato \cite{s} that $\int_{-\infty}^{\infty}|y|\xi_2(dy)< \infty$.
Thus, by dominated convergence theorem,
$$I_2 \sim p_1(x)\int_{-\infty}^{0}|y|\xi_2(dy).$$
 Note from Lemma 4.1 that $p_1(x)\sim p(x)$. Hence we see that
\begin{equation}
\begin{split}
\bar\mu(x)&=\bar\xi_1(x)+m(\xi_2)p_1(x) +o(p_1(x))\\
&=\bar\xi_1(x)+m(\xi_2)p(x) +o(p(x))
\end{split}
\end{equation}
as $x \to \infty$. Note that $p^{2\otimes}(x)\sim 2 p(x)$ and $m(\xi_2^{2*})=2m(\xi_2).$ In the same way, we have
\begin{equation}
\begin{split}
\overline{\mu^{2*}}(x)&=\overline{\xi_1^{2*}}(x)+m(\xi_2^{2*})p^{2\otimes}(x) +o(p^{2\otimes}(x))\\
&=\overline{\xi_1^{2*}}(x)+4m(\xi_2)p(x) +o(p(x))\nonumber
\end{split}
\end{equation}
as $x \to \infty$. Hence we obtain from (4.7) that
\begin{equation}
\begin{split}
&\overline{\mu^{2*}}(x)-2\overline{\mu}(x)\\
&=\overline{\xi_1^{2*}}(x)-2\overline{\xi_1}(x)+2m(\xi_2)p(x) +o(p(x)).\nonumber
\end{split}
\end{equation}
as $x \to \infty$. Since $\int_{-\infty}^{\infty}|y|\xi_2(dy)< \infty$, we see from Theorem 25.3 of Sato \cite{s} that $\int_{-\infty}^{\infty}|x|\mu(dx) < \infty$ if and only if $0< m(\xi_1) < \infty.$ Thus we have
$$\overline{\mu^{2*}}(x)=2\overline{\mu}(x)+2m(\mu)p(x) +o(p(x))$$
 as $x \to \infty$
 if and only if
$$\overline{\xi_1^{2*}}(x)=2\overline{\xi_1}(x)+2m(\xi_1)p_1(x)+
o(p_1(x))$$
as $x \to \infty$. Thus under the assumption of $\int_{-\infty}^{0}|y|\mu(dy)< \infty$, we have
$p(x) \in \mathcal{S}^2_{d}$ if and only if $p_1(x) \in \mathcal{S}^2_{d}$, equivalently 
$\xi_1 \in \mathcal{S}^2_{loc}$. We find from Theorem 1.1 that
$\xi_1 \in \mathcal{S}^2_{loc}$ if and only if $\nu_{(1)}  \in \mathcal{S}^2_{loc}$. That is, $p(x) \in \mathcal{S}^2_{d}$ if and only if $\nu_{(1)}  \in \mathcal{S}^2_{loc}$, equivalently, $\frac{1}{\bar\nu(1)}1_{(1,\infty)}(x)k(x)/x \in \mathcal{S}^2_{d}$. Next we prove assertion (ii). If $p(x) \in \mathcal{S}^2_{d}$, then  $\nu_{(1)}  \in \mathcal{S}^2_{loc}$ and hence by Theorem 1.1 we have
\begin{equation}
\begin{split}
\bar\xi_1(x)&=\bar\nu(x)+m(\xi_1)\nu((x,x+1])+o(\nu((x,x+1])\\ 
&=\bar\nu(x)+m(\xi_1)p(x)+o(p(x))\\
&=\bar\nu(x)+m(\xi_1)p_1(x)+o(p_1(x))
\end{split}
\end{equation}
as $x \to \infty.$ Thus it follows from (4.7) that (4.5) and (4.6) hold. Next we prove assertion (iii). The assumption that (4.5) with finite $m(\mu)$, $p(x) \in \mathcal{S}^2_{d} $, and $(\bar\mu(x))^2=o(p(x))$ as $x \to \infty$ implies that
(4.8) with finite $m(\xi_1)$, $\xi_1 \in \mathcal{S}_{loc}$, and $(\bar\xi_1(x))^2=o(\xi_1((x,x+1]))$ as $x \to \infty$. Thus we see from (iii) of Theorem 1.1 that 
 $\xi_1 \in \mathcal{S}^2_{ loc}$, that is, $p_1(x) \in \mathcal{S}^2_{ d}$. It follows from the proof of (i) that
$p(x) \in \mathcal{S}^2_{d}$. \hfill $\Box$

\begin{cor}
Let  $\mu(dx)=p(x)dx$  be a self-decomposable  distribution on $\mathbb R$  with L\'evy measure $\nu$.   Then, the following hold : 

(i) $p(x) \in \mathcal{S}^2_{d}$ if and only if $p^{t}(x) \in \mathcal{S}^2_{d}$ for some $t>0$, equivalently, for all $t>0$.

(ii) If $p(x) \in \mathcal{S}^2_{d}$, then, for all $t>0$,
\begin{equation}
\overline{\mu^{t*}}(x)=t\bar\mu(x)+(t^2-t)m(\mu)p(x)+o(p(x))\nonumber
\end{equation} 
as $x \to \infty.$ 
\end{cor}
Proof. By argument analogous to the proof of Corollary 1.1, we can easily prove the corollary from Theorem 4.1 and Lemma 4.2.  \hfill $\Box$
\section{Examples}
By using a method of  Kl\"uppelberg \cite{k1} and Baltrunas \cite{b},
Lin \cite{l} proved that the standard lognormal distribution,  Weibull distribution with parameter $\beta \in(0,1)$, and Pareto distribution with parameter $\alpha >1$ belong to the class $\mathcal{S}^2_{loc}$. Those distributions are all self-decomposable, so their densities also belong to the class $\mathcal{S}^2_{d}$.
See Sato \cite{s} and Steutel and van Harn \cite{sh} for their self-decomposability. The following examples are direct consequence of Theorem 1.1
and Corollary 1.1 and hence  their
 proofs are omitted.
\begin{exm} Let $\mu$  be the standard lognormal distribution with L\'evy measure $\nu(dx)=k(x)/xdx$. Then we have the density 
$$p(x):=\frac{1}{\sqrt{2\pi}x}\exp\left(-\frac{(\log x)^2}{2}\right)$$
for $x >0$. Embrechts et al. \cite{egv} showed that $\mu$ is subexponential and that
$$\bar\nu(x)\sim \bar\mu(x)\sim \frac{x}{\log x}p(x)$$
and
$$\overline{\mu^{t*}}(x)\sim t\bar\mu(x).$$
Watanabe and Yamamuro \cite {wy2} proved a conjecture of Bondesson \cite{bo}. That is, 
$$k(x)\sim xp(x).$$
We have
\begin{equation}
\bar\nu(x)=\bar\mu(x)\left(1-\sqrt{e}\frac{\log x}{x} +o\left(\frac{\log x}{x}\right)\right) \nonumber
\end{equation} 
as $x \to \infty$ and, for $t >0$,
\begin{equation}
\overline{\mu^{t*}}(x)=t\bar\mu(x)\left(1+(t-1)\sqrt{e}\frac{\log x}{x} +o\left(\frac{\log x}{x}\right)\right)\nonumber
\end{equation} 
as $x \to \infty.$
\end{exm}
\begin{exm} Let $\mu$ be  Weibull distribution with L\'evy measure $\nu$ and parameter $\beta \in(0,1)$. Then we have 
$$\bar\mu(x):=\exp(-x^{\beta})$$
 for $ x \in \mathbb R_+$,
\begin{equation}
\bar\nu(x)=\bar\mu(x)(1-\Gamma(\beta^{-1})x^{\beta-1} +o(x^{\beta-1}))\nonumber
\end{equation} 
as $x \to \infty$, and, for $t >0$,
\begin{equation}
\overline{\mu^{t*}}(x)=t\bar\mu(x)(1+(t-1)\Gamma(\beta^{-1})x^{\beta-1} +o(x^{\beta-1}))\nonumber
\end{equation} 
as $x \to \infty.$
\end{exm}
\begin{exm} Let $\mu$ be  Pareto distribution with L\'evy measure $\nu$ and parameter $\alpha >1$. Then we have 
$$\bar\mu(x):=(1+x)^{-\alpha}$$
 for $x \in \mathbb R_+$,
\begin{equation}
\bar\nu(x)=\bar\mu(x)\left(1-\frac{\alpha}{\alpha-1}x^{-1} +o(x^{-1})\right)\nonumber
\end{equation} 
as $x \to \infty,$ and, for $t >0$,
\begin{equation}
\overline{\mu^{t*}}(x)=t\bar\mu(x)\left(1+(t-1)\frac{\alpha}{\alpha-1}x^{-1} +o(x^{-1})\right)\nonumber
\end{equation} 
as $x \to \infty.$
\end{exm}
\section{Remarks on the regularly varying case}
 We cannot find from our results the relations of Example 5.3 for Pareto distribution with parameter $0< \alpha \leq 1$ because it does not belong to the class $\mathcal{S}^2_{loc}$. 
However, we can get the analogous relations by using the following lemma of Omey and Willekens \cite{ow}. Theorem 4.3 of \cite{ow} is a direct consequence from Theorem 2.3 of \cite{ow} for a compound Poisson distribution on $\mathbb R_+$, but there is a mistake in the case of finite mean for 
 an  infinitely divisible distribution on $\mathbb R_+$. So we restore and prove it for 
 an  infinitely divisible distribution on $\mathbb R_+$.
\begin{lem}(Theorem 4.3 of \cite{ow})
Let $\mu$ be an  infinitely divisible distribution on $\mathbb R_+$  with L\'evy measure $\nu$. Assume
 that $\nu(dx)$ has a density $q(x)$ on $(1, \infty)$ such that
 $q(x)\sim  x^{-\alpha-1}l(x)$ for $0 \leq \alpha \leq 1$ with $l(x)$ being slowly varying as $x \to \infty$. Define a constant $C(\alpha)$ for $0< \alpha < 1$ as 
$$C(\alpha):=\frac{(1-\alpha)(2\alpha-1)(\Gamma(1-\alpha))^2}{2\alpha\Gamma(2-2\alpha)}.$$

(i) We have for $0< \alpha < 1$ 
\begin{equation}
\lim_{x \to \infty}\frac{\bar\mu(x)-\bar\nu(x)}{q(x)\int_1^x\bar\nu(u)du}=C(\alpha).
\end{equation} 

(ii) For $ \alpha = 1$, if $\int_1^{\infty}\bar\nu(u)du= \infty$, then we have
\begin{equation}
\lim_{x \to \infty}\frac{\bar\mu(x)-\bar\nu(x)}{q(x)\int_1^{x}\bar\nu(u)du}=1.
\end{equation} 

(iii) For $ \alpha = 1$, if $\int_1^{\infty}\bar\nu(u)du< \infty$, 
then we have
\begin{equation}
\lim_{x \to \infty}\frac{\bar\mu(x)-\bar\nu(x)}{q(x)m(\mu)}=1.
\end{equation} 

(iv) For $ \alpha = 0$, 
then we have
\begin{equation}
\lim_{x \to \infty}\frac{\bar\mu(x)-\bar\nu(x)}{(\bar\nu(x))^2}=-\frac{1}{2}.
\end{equation}
\end{lem}
Proof.  Let $\mu$ be an  infinitely divisible distribution on $\mathbb R_+$  with L\'evy measure $\nu$. Assume
 that $\nu(dx)$ has a density $q(x)$ on $(1, \infty)$ such that
 $q(x)\sim  x^{-\alpha-1}l(x)$ for $0\leq \alpha \leq 1$ with $l(x)$ being slowly varying as $x \to \infty$. Define a  compound Poisson distribution $\mu_1$ on $\mathbb R_+$ as (2.2) for $c=1$.
 Define
an infinitely divisible distributions $\mu_2$ on $\mathbb R_+$ as $\mu=\mu_1*\mu_2$. Then we have by Theorem 2.3 of  \cite{ow}, for $0\leq \alpha \leq 1$, the lemma is true by substituting $\mu_1$ for $\mu$. Thus we can assume that $\mu_2(dx)\ne \delta_0(dx)$.
We see from Theorem 25.17  of Sato \cite{s} that, for every $b>0$,  $\int_{0-}^{\infty}\exp(b x)\mu_2(dx) < \infty$ and hence $\overline{\mu_2}(x) = o(e^{-bx})$ as $x \to \infty$. 
We have
\begin{equation}
\begin{split}
&\bar\mu(x)-\bar\mu_1(x)\\
&=\overline{\mu_1*\mu_2}(x)-\bar\mu_1(x)\\
&=\int_{0-}^{\infty}\mu_1((x-y,x])\mu_2(dy)\\
&=I_1+I_2 +I_3,\nonumber
\end{split}
\end{equation}
where 
$$I_1:=\int_{0-}^{A+}\mu_1((x-y,x])\mu_2(dy),$$
$$I_2:=\int_{A+}^{x/2+}\mu_1((x-y,x])\mu_2(dy),$$
and
$$I_3:=\int_{x/2+}^{\infty}\mu_1((x-y,x])\mu_2(dy).$$
Since $q(x)\sim  x^{-\alpha-1}l(x)$, $\nu_{(1)}\in \mathcal{S}_{ loc}$ and hence, by Lemma 2.1, $\mu_1\in \mathcal{S}_{ loc}$. Thus,
\begin{equation}
\begin{split}
I_1 &\sim\mu_1((x,x+1])\int_{0-}^{A+}y\mu_2(dy)\\  \nonumber
& \sim \mu_1((x,x+1])\int_{0-}^{\infty}y\mu_2(dy)
\end{split}
\end{equation}
as $x \to \infty$ and then $A \to \infty$. 
Since $\mu_1\in \mathcal{S}_{ loc}$, there are
$C > 0$
and $\epsilon >0$ such that, for $0 \leq y \leq x/2$ and for sufficiently large $x >0$,
$$\mu_1((x-y,x])\leq Ce^{\epsilon y}\mu_1((x,x+1]).$$
 Thus we have
\begin{equation}
\begin{split}
I_2 &\leq \mu_1((x,x+1])\int_{A+}^{x/2+}Ce^{\epsilon y}\mu_2(dy)\\
&=o(\mu_1((x,x+1]))\nonumber
\end{split}
\end{equation}
as $x \to \infty$ and then $A \to \infty$. 
$$I_3\leq\bar\mu_2(x/2)=o(e^{-x})=o(\mu_1((x,x+1]))$$
as $x \to \infty$. Thus we see that
\begin{equation}
\bar\mu(x)-\bar\mu_1(x) \sim m(\mu_2)\mu_1((x,x+1]).
\end{equation}
Note from Lemma 4.1 that, for $0 < \alpha <1$ or $\alpha =1$ with $\int_1^{\infty}\bar\nu(u)du=\infty$, 
$$\mu_1((x,x+1])\sim q(x)=o(q(x)\int_1^x\bar\nu(u)du)$$
as $x \to \infty$. For $\alpha =0$, we have by Lemma 4.1
$$\mu_1((x,x+1])\sim q(x)=o((\bar\nu(x))^2)$$
as $x \to \infty$. Thus except the case of $\alpha =1$ with finite $m(\mu_1)$, the lemma is true. In the case of $\alpha =1$ with finite $m(\mu_1)$, we see from  (6.3) with   substituting $\mu_1$ for $\mu$ and (6.5) that the lemma is true.\hfill $\Box$

\begin{prop}
Let $\mu(dx)=p(x)dx$ be a self-decomposable distribution on $\mathbb R_+$. Assume that $p(x) \sim x^{-\alpha-1}l(x)$ for $0\leq\alpha \leq 1$ with $l(x)$ being slowly varying as $x \to \infty$. Define  slowly varying functions $l^*(x)$ and $l_*(x)$
as $l^*(x)=\int_1^{x}l(u)/u du$ and $l_*(x)=\int_x^{\infty}l(u)/u du$ for $x >1$. Then we have the following :

(i) Let $0< \alpha < 1$ and define $K(\alpha)$ as
$$K(\alpha):=\frac{(2\alpha-1)(\Gamma(1-\alpha))^2}{2\alpha\Gamma(2-2\alpha)}.$$
 Then we have
\begin{equation}
\bar\nu(x)=\bar\mu(x)\left(1-K(\alpha)x^{-\alpha}l(x) +o(x^{-\alpha}l(x))\right)
\end{equation} 
as $x \to \infty,$ and, for $t >0$,
\begin{equation}
\overline{\mu^{t*}}(x)=t\bar\mu(x)\left(1+(t-1)K(\alpha)x^{-\alpha}l(x) +o(x^{-\alpha}l(x))\right)
\end{equation} 
as $x \to \infty.$

(ii) Let $ \alpha = 1$. Assume that $l^*(\infty)=\infty$. Then we have
\begin{equation}
\bar\nu(x)=\bar\mu(x)\left(1-\frac{l^*(x)}{x} +o\left(\frac{l^*(x)}{x}\right)\right)
\end{equation} 
as $x \to \infty,$ and, for $t >0$,
\begin{equation}
\overline{\mu^{t*}}(x)=t\bar\mu(x)\left(1+(t-1)\frac{l^*(x)}{x} +o\left(\frac{l^*(x)}{x}\right)\right)
\end{equation} 
as $x \to \infty.$

(iii) Let $ \alpha = 1$. Assume that $l^*(\infty)<\infty$. Then we have
\begin{equation}
\bar\nu(x)=\bar\mu(x)\left(1-\frac{m(\mu)}{x} +o\left(\frac{1}{x}\right)\right)
\end{equation} 
as $x \to \infty,$ and, for $t >0$,
\begin{equation}
\overline{\mu^{t*}}(x)=t\bar\mu(x)\left(1+(t-1)\frac{m(\mu)}{x} +o\left(\frac{1}{x}\right)\right)
\end{equation} 
as $x \to \infty.$

(iv)
Let $ \alpha =0$.
Then we have
\begin{equation}
\bar\nu(x)=\bar\mu(x)\left(1+\frac{l_*(x)}{2} +o(l_*(x))\right)
\end{equation} 
as $x \to \infty,$ and, for $t >0$,
\begin{equation}
\overline{\mu^{t*}}(x)=t\bar\mu(x)\left(1-(t-1)\frac{l_*(x)}{2}  +o(l_*(x))\right)
\end{equation}
\end{prop} 
Proof.  Assume that $p(x) \sim x^{-\alpha-1}l
(x)$ for $0\leq\alpha \leq 1$ with $l(x)$ being slowly varying as $x \to \infty$.  
First we prove (i). Let $0< \alpha < 1$. Since
$p(x)\in \mathcal{S}_{d}$, we have by Lemma 4.1
$$q(x)\sim  x^{-\alpha-1}l(x).$$ 
By Karamata's theorem(Theorem 1.5.11  of \cite{bgt}), we have
$$\bar\nu(x)\sim \bar\mu(x)\sim \frac{x^{-\alpha}l(x)}{\alpha}$$
and
$$\int_1^x\bar\nu(u)du \sim \frac{x^{1-\alpha}l(x)}{
\alpha(1-\alpha)}.$$
Thus we see from  (6.1) of Lemma 6.1 that
$$\lim_{x \to \infty}\frac{\bar\mu(x)-\bar\nu(x)}{x^{-2\alpha}(l(x))^2}=\frac{K(\alpha)}{\alpha}
.$$
Thus we have (6.6). In the same way, we have
$$\lim_{x \to \infty}\frac{\overline{\mu^{t*}}(x)-t\bar\nu(x)}{x^{-2\alpha}(l(x))^2}=t^2\frac{K(\alpha)}{\alpha}.$$
Hence we get (6.7) by (6.6).
Next we prove (ii).   Assume that $p(x) \sim x^{-2}
l(x)$. Then, by Karamata's theorem, we have $\bar\mu(x) \sim x^{-1}l(x)$. We have by Lemma
4.1
$$q(x)\sim  x^{-2}l(x).$$
We see from Karamata's theorem that $\bar\nu(x) \sim x^{-1}l(x)$ and
$$\int_1^x\bar\nu(u)du \sim l^*(x)$$
and that $\int_1^{\infty}\bar\nu(u)du= \infty$ from $l^*(\infty)=\infty$. Thus we see from (6.2) of Lemma 6.1 that
$$\lim_{x \to \infty}\frac{\bar\mu(x)-\bar\nu(x)}{x^{-2}l(x)l^*(x)}=1.$$
Thus we have (6.8). In the same way, we have
$$\lim_{x \to \infty}\frac{\overline{\mu^{t*}}(x)-t\bar\nu(x)}{x^{-2}l(x)l^*(x)}=t^2.$$
Hence we get (6.9)
 by (6.8).
 Next we prove (iii). As in (ii), we have $q(x)\sim p(x) \sim x^{-2}
l(x)$, $\bar\nu(x) \sim \bar\mu(x) \sim x^{-1}l(x)$, and 
$$\int_1^x\bar\nu(u)du \sim l^*(x).$$
 We see that
 $\int_1^{\infty}\bar\nu(u)du< \infty$ from $l^*(\infty)<\infty$. Thus we find from (6.3) of Lemma 6.1 that
$$\lim_{x \to \infty}\frac{\bar\mu(x)-\bar\nu(x)}{x^{-2}l(x)m(\mu)}=1.$$
Thus we have (6.10). In the same way, we have
$$\lim_{x \to \infty}\frac{\overline{\mu^{t*}}(x)-t\bar\nu(x)}{x^{-2}l(x)m(\mu)}=t^2.$$
Hence we get (6.11) by (6.10). Next we prove (iv).  Assume that $p(x) \sim x^{-1}
l(x)$. Then, we see from Lemma 4.1
that $q(x) \sim x^{-1}
l(x)$. Thus we have
$$\bar\mu(x)\sim \bar\nu(x)\sim l_*(x).$$
We find from (6.4) of Lemma 6.1 that 
$$\lim_{x \to \infty}\frac{\bar\mu(x)-\bar\nu(x)}{(l_*(x))^2}=-\frac{1}{2}.$$
Thus we have (6.12). In the same way, we have
$$\lim_{x \to \infty}\frac{\overline{\mu^{t*}}(x)-t\bar\nu(x)}{(l_*(x))^2}=-\frac{t^2}{2}.$$
Hence we get (6.13) by (6.12).
\hfill $\Box$

Finally, we give the relations for Pareto distribution with parameter $0< \alpha \leq 1$ as an example of Proposition 6.1. They are different from the relations of Example 5.3.
\begin{exm} Let $\mu$ be  Pareto distribution with L\'evy measure $\nu$ and parameter $0< \alpha \leq 1$. Then we have 
$$\bar\mu(x):=(1+x)^{-\alpha}$$
 for $x \in \mathbb R_+$.
 
(i) Let $0< \alpha < 1$.
Then we have
\begin{equation}
\bar\nu(x)=\bar\mu(x)\left(1-\alpha K(\alpha)x^{-\alpha} +o(x^{-\alpha})\right)\nonumber
\end{equation} 
as $x \to \infty,$ and, for $t >0$,
\begin{equation}
\overline{\mu^{t*}}(x)=t\bar\mu(x)\left(1+(t-1)\alpha K(\alpha)x^{-\alpha} +o(x^{-\alpha})\right)\nonumber
\end{equation} 
as $x \to \infty.$

(ii) Let $ \alpha = 1$.  Then we have
\begin{equation}
\bar\nu(x)=\bar\mu(x)\left(1-\frac{\log x}{x} +o\left(\frac{\log x}{x}\right)\right)\nonumber
\end{equation} 
as $x \to \infty,$ and, for $t >0$,
\begin{equation}
\overline{\mu^{t*}}(x)=t\bar\mu(x)\left(1+(t-1)\frac{\log x}{x} +o\left(\frac{\log x}{x}\right)\right)\nonumber
\end{equation} 
as $x \to \infty.$
\end{exm}

 Toshiro Watanabe\\
 Center for Mathematical Sciences, The University of Aizu, Aizu-Wakamatsu Fukushima 965-8580, Japan \\
e-mail: t-watanb@u-aizu.ac.jp

\end{document}